\documentclass[letterpaper,12pt]{article}

\usepackage[letterpaper,top=3cm,bottom=3cm,left=3cm,right=3cm,marginparwidth=1.75cm]{geometry}

\usepackage[english]{babel}
\usepackage[utf8x]{inputenc}
\usepackage[T1]{fontenc}

\usepackage{amsmath}
\usepackage{amsthm}
\usepackage{amssymb}
\usepackage[all,cmtip]{xy}
\usepackage{graphicx}
\usepackage{theoremref}
\usepackage{hyperref}

\theoremstyle{plain}
\newtheorem{theorem}{Theorem}[section]
\newtheorem{lemma}[theorem]{Lemma}
\newtheorem{prop}[theorem]{Proposition}
\newtheorem{cor}[theorem]{Corollary}

\theoremstyle{definition}
\newtheorem{exx}[theorem]{Example}

\usepackage{theoremref}

\newcommand{\Z}{\mathbb{Z}}
\newcommand{\Q}{\mathbb{Q}}

\newcommand{\K}{\mathbb{K}}

\newcommand{\BAG}{B^A(G)}
\newcommand{\BAH}{B^A(H)}
\newcommand{\gBAG}{\widetilde{B^A}(G)}
\newcommand{\BAKG}{B^A_{\mathbb{K}}(G)}

\newcommand{\BAOG}{B^A_{\mathcal{O}}(G)}
\newcommand{\gBAOG}{\widetilde{B^A_{\mathcal{O}}}(G)}

\newcommand{\MAG}{\mathcal{M}^A_G}

\newcommand{\DAG}{\mathcal{D}^A_G}

\newcommand{\DnG}{\mathcal{D}^{C_n}_G}

\newcommand{\SG}{\mathcal{S}_G}

\newcommand{\eGAGPhi}{e^{G,A}_{G,\Phi}}

\newcommand{\sGAHPhi}{s^{G,A}_{H,\Phi}}
\newcommand{\eGAHPhi}{e^{G,A}_{H,\Phi}}
\newcommand{\cGAHPhi}{c^{G,A}_{H,\Phi}}

\newcommand{\eGnHh}{e^{G,n}_{H,h}}
\newcommand{\cGnHh}{c^{G,n}_{H,h}}

\title{The conductors of fibered Burnside rings}

\author{
Benjam\'in Garc\'ia\\
Divisi\'on de Carreras de Ingenier\'ias\\
Universidad Tecnol\'ogica de San Miguel de Allende\\
San Miguel de Allende, Guanajuato 37700\\
M\'exico\\
\texttt{bagh1704@gmail.com}
}

\date{\today}

\begin{document}

\maketitle

\begin{abstract}
We give an explicit formula for the conductors of the primitive idempotents of fibered Burnside algebras of finite rank. This formula generalizes previous results for Burnside rings, monomial representation rings, and cyclic fibers, removing the coprimality condition required in the latter case. We also derive some consequences for cyclic fiber groups and $p$-groups.
\end{abstract}



\section*{Introduction}

The conductors of the primitive idempotents of the rational Burnside algebra were first determined by Nicolson in 1978 \cite{Nicolson}. In his article, Nicolson showed that, for a primitive idempotent $e^{B(G)}_H$ of $B_{\Q}(G)$ (indexed by the conjugacy class of $H\leq G$), the least positive integer $t$ such that $te^{B(G)}_H$ belongs to $B(G)$ is precisely $[N_G(H):H][H:H']_0$, where $[H:H']_0$ denotes the product of the prime divisors of $[H:H']$. Later, Kratzer \& Thévenaz \cite{KraThe} and Dress \& Vallejo \cite{DressVall} provided more direct proofs of this result. 

In the case of the monomial representation algebra $D_{\K}(G)$, the primitive idempotents $e^{D(G)}_{H,h}$ are indexed by the conjugacy classes of pairs $(H,hH')$ with $H\leq G$ and $h\in H$, where $\K$ is a field of characteristic zero containing a primitive $\mathrm{exp}(G)$-th root of unity. Müller showed in \cite[\S 3]{Muller} that the conductor of $e^{D(G)}_{H,h}$ is equal to the index $\left[N_G\left(H,hH'\right):H'\right]$. In general, $\left[N_G\left(H,hH'\right):H'\right]e^{D(G)}_{H,h}$ does not lie in $D(G)$ but in $D_{\mathcal{O}}(G)$, where $\mathcal{O}$ is the ring of integers of $\K$.

A first generalization of the previous results in the context of fibered Burnside rings was given by Raggi-Cárdenas and the author in \cite[\S 5]{GarRag}. For a positive integer $n$, the primitive idempotents $e^{G,n}_{H,h}$ of $B^{C_n}_{\K}(G)$, where $\K$ is a field of characteristic zero containing a primitive $n$-th root of $1$, are indexed by the conjugacy classes of pairs $\left(H,hH^{[n]}\right)$, with $H\leq G$ and $h\in H$. Here, $H^{[n]}$ is the minimal normal subgroup of $H$ such that $H/H^{[n]}$ is an abelian group of exponent dividing $n$ \cite[\S 1]{GarRag}. When $n$ is a divisor of $\mathrm{exp}(G)$ such that $\gcd(n,\mathrm{exp}(G)/n)=1$, the conductor of $\eGnHh$ is precisely $\left[N_G\left(H,hH^{[n]}\right):H^{[n]}\right]\left[H^{[n]}:H'\right]_0$. Although this result has restrictive hypotheses, it covers the cases of the Burnside ring ($n=1$) and the monomial representation ring ($n=\mathrm{exp}(G)$).

In this article, we focus on fibered Burnside rings $\BAG$ of finite rank. As shown by Boltje and Y\i lmaz in \cite[\S 3]{BolYil}, the primitive idempotents $\eGAHPhi$ of the fibered Burnside algebra $\BAKG$ over a sufficiently large field of characteristic zero are indexed by orbits of pairs $(H,\Phi)$ under a natural action of $G$ on the set $\DAG$. The main result provides an explicit formula for the conductor of $\eGAHPhi$, namely,
$$\cGAHPhi=\left[N_G(H,\Phi):H\right]\left|\mathrm{Hom}(H,A)\right|r^A(H),$$ 
where $r^A(H)$ is an additional factor depending only on $H$ and $A$. As in the previous cases, Boltje's integrality criterion and the Möbius function play an important role in the proof. A closer examination of the values of the Möbius function allows us to remove the coprimality condition imposed in \cite{GarRag}.

The article is organized into two sections. The first collects the necessary preliminaries on abelian groups, the Möbius function, and fibered Burnside rings. The second is devoted to the statement and proof of the main result, together with some consequences concerning cyclic fibers and $p$-groups.

\subsection*{Notation}

For a positive integer $n$, we write $n_0$ for its largest square-free divisor, that is, the product of its distinct prime divisors. We write $\mathcal{S}_G$ for the poset of subgroups of a group $G$, and $[\mathcal{S}_G]$ for a subset of representatives of the conjugacy classes of $\SG$. If $X$ is a $G$-set, then $X/G$ denotes its set of orbits, and we write $x=_G y$ for elements $x,y\in X$ lying in the same orbit. We regard $G$ and $\mathcal{S}_G$ as $G$-sets under conjugation, writing ${}^g x$ instead of $gxg^{-1}$. For $H\leq G$ and an abelian group $A$, $\phi\in\mathrm{Hom}(H,A)$, and $g\in G$, we write ${}^g\phi$ for the morphism defined by $({}^g\phi)(x)=\phi(g^{-1}xg)$, for $x\in {}^gH$. For $K\leq H$, we denote by $\phi|_K$ the restriction of $\phi$ to $K$. All rings $R$ are assumed to be associative, unital, and commutative, and we use the standard notation $\mathbb{Z}$, $\mathbb{Q}$, $\mathbb{R}$, and $\mathbb{C}$ for the usual number systems.

\section{Preliminaries}

We recall some facts about morphisms to abelian groups and the Möbius function of a poset of subgroups, and then review the basic properties of fibered Burnside rings that will be used throughout the paper.

\subsection{Morphisms to abelian groups}

For an abelian group $A$ and a positive integer $n$, we denote by $\mathrm{tor}_n(A)$ the subgroup of elements of $A$ whose orders divide $n$; we say that $A$ is an \textit{$n$-torsion} abelian group if $\mathrm{tor}_n(A)=A$. For a finite group $G$, the subgroup $G^{[n]}$ is defined as the smallest normal subgroup $T\unlhd G$ such that $G/T$ is an $n$-torsion abelian group; it can be constructed as the subgroup generated by $G'\cup \{g^n \mid g\in G\}$. By \cite[Lemma 1.1]{GarRag}, $G^{[\mathrm{exp}(G)]}=G'$, $G^{[1]}=G$, $H^{[n]}\leq G^{[n]}$ for $H\leq G$, and $G^{[n]}\leq G^{[m]}$ if $m\mid n$. Also, there is a natural isomorphism 
$$\mathrm{Hom}\left(G,\mathrm{tor}_n(A)\right) \cong \mathrm{Hom}\left(G/G^{[n]},A\right)$$
by \cite[Lemma 1.2]{GarRag}. In particular,
$$\mathrm{Hom}\left(G,A\right)\cong \mathrm{Hom}\left(G,\mathrm{tor}_{\mathrm{exp}(G)}(A)\right) \cong\mathrm{Hom}\left(G/G^{[n(A,G)]},A\right),$$
where $n(A,G):=\mathrm{exp}\left(\mathrm{tor}_{\mathrm{exp}(G)}(A)\right)$. By \cite[Lemma 1.3]{GarRag}, $\mathrm{Hom}(G,C_n)$ and $G/G^{[n]}$ are isomorphic, and so $\left|\mathrm{Hom}(G,C_n)\right|=\left[G:G^{[n]}\right]$. 

\begin{lemma}\thlabel{discardprimes}
    Let $n$ be a positive integer and $G$ be a finite abelian group. Let $p$ be a prime divisor of $\left|G^{[n]}\right|$.
    \begin{enumerate}
        \item If $p$ divides $\left[G:G^{[n]}\right]$, then $G$ has a $p$-element of order greater than $p$;
        \item if $p$ does not divide $\left[G:G^{[n]}\right]$, then $p$ does not divide $n$.
    \end{enumerate}
\end{lemma}
\begin{proof}
    Since $G$ is finite and abelian, the map $f_n:G\to G$ given by $g\mapsto g^n$ is an endomorphism with $\mathrm{Im}(f_n)=G^{[n]}$, and $\left|\mathrm{ker}(f_n)\right|=\left[G:G^{[n]}\right]$. Because $p\mid \left|G^{[n]}\right|$, the $p$-primary component $G_p$ is non-trivial, and $f_n$ restricts to an endomorphism on $G_p$ with image $G^{[n]}_p=G^{[n]}\cap G_p\neq 1$. If $p\mid \left[G:G^{[n]}\right]$, then $p$ divides $\mathrm{exp}\left(G/G^{[n]}\right)$, which divides $n$. Then there is some $g\in G_p$ such that $1\neq g^n\in G^{[n]}_p$; since $p\mid n$, it follows that $g$ has order greater than $p$. On the other hand, if $p\nmid \left[G:G^{[n]}\right]$, then $G_p\cap \mathrm{ker}(f_n)=1$, hence $f_n$ is an automorphism on $G_p$ and so $n$ is coprime to $p$.
\end{proof}

If $A$ is a finite abelian group, we consider a decomposition $A=C_{d_1}\times C_{d_2}\times \cdots \times C_{d_r}$ into cyclic factors. For any finite group $G$, there is a natural isomorphism $\text{Hom}(G,A)\cong \prod_{i=1}^r\text{Hom}(G,C_{d_i})$, and for each $H\leq G$, the following diagram involving the restriction maps commutes:
    $$\xymatrix{
        \mathrm{Hom}(G,A)\ar[d]_{\mathrm{res}^G_H}\ar[rr] &&\prod_{i=1}^r\mathrm{Hom}(G,C_{d_i})\ar[d]^{\prod_i\text{res}^G_H}\\
        \mathrm{Hom}(H,A)\ar[rr] &&\prod_{i=1}^r\mathrm{Hom}(H,C_{d_i}).
    }$$
Indeed, for every $\phi:G\longrightarrow A$ and all $1\leq i \leq r$, $\pi_{C_{d_i}}\circ (\phi|_H)=(\pi_{C_{d_i}}\circ \phi)|_H$.

\begin{lemma}\thlabel{MullerLemma}
    Let $H\leq G$, and consider the restriction map $\mathrm{res}^G_H:\mathrm{Hom}(G,A)\longrightarrow \mathrm{Hom}(H,A)$. The following holds:
    \begin{enumerate}
        \item For $\psi\in \mathrm{Im}\left(\mathrm{res}^G_H\right)$, 
        $$\left|\left(\mathrm{res}^G_H\right)^{-1}(\psi)\right|=\prod_{i=1}^r\left[G:HG^{[d_i]}\right];$$
        \item Let $\Phi:\mathrm{Hom}(G,A)\longrightarrow \K^{\times}$ be a linear $\K$-character, where $\K$ is a field of characteristic zero. Then
        $$\sum_{\phi\in \left(\mathrm{res}^G_H\right)^{-1}(\psi)}\Phi(\phi)=\begin{cases}
            \left(\prod_{i=1}^r\left[G:HG^{[d_i]}\right]\right)\Phi(\phi_0) &\text{if }\mathrm{ker}\left(\mathrm{res}^G_H\right)\leq \mathrm{ker}\;\Phi,\\
            0 &\text{otherwise,}
        \end{cases}$$
        where $\phi_0$ is any element in $\left(\mathrm{res}^G_H\right)^{-1}(\psi)$.
    \end{enumerate}
\end{lemma}
\begin{proof}
    Note that $\left(\mathrm{res}^G_H\right)^{-1}(\psi)=\mathrm{ker}\left(\mathrm{res}^G_H\right)\phi_0$ for some $\phi_0:G\longrightarrow A$ such that $\phi_0|_H=\psi$, and that
    $$\mathrm{ker}\left(\mathrm{Hom}(G,A)\to \mathrm{Hom}(H,A)\right)\cong \prod_{i=1}^r\ker\left(\mathrm{Hom}(G,C_{d_i})\to \mathrm{Hom}(H,C_{d_i})\right)$$
    by the previous diagram. By \cite[Lemma 5.3]{GarRag}, 
    $$\left|\mathrm{ker}\left(\mathrm{Hom}\left(G,C_{d_i}\right)\to \mathrm{Hom}\left(H,C_{d_i}\right)\right)\right|=\left[G:HG^{[d_i]}\right]$$
    for each $i$, proving assertion 1. 
    
    For the second assertion, we have
    \begin{align*}
        \sum_{\phi\in \left(\mathrm{res}^G_H\right)^{-1}(\psi)}\Phi(\phi) 
        &=\sum_{\lambda \in \mathrm{ker}\left(\mathrm{res}^G_H\right)}\Phi(\lambda\phi_0)\\
        &=\left(\sum_{\lambda\in \mathrm{ker}\left(\mathrm{res}^G_H\right)}\Phi(\lambda)\right)\Phi(\phi_0)\\
        &=\begin{cases}
            \left|\mathrm{ker}\left(\mathrm{res}^G_H\right)\right|\Phi(\phi_0) &\text{if }\mathrm{ker}\left(\mathrm{res}^G_H\right)\leq \mathrm{ker}\;\Phi,\\
            0 &\text{otherwise,}
        \end{cases}
    \end{align*}
    by orthogonality relations for $\mathrm{ker}\left(\mathrm{res}^G_H\right)$. Together with assertion 1, this proves 2.
\end{proof}

\subsection{The Möbius function of the poset of subgroups}

We consider $\SG$ as a poset with respect to inclusion. The \textit{Möbius function} $\mu:\SG\times \SG\longrightarrow \Z$ is defined recursively as follows: for $H\leq G$, $\mu(H,H)=1$, and for $K\neq H$,
\[\mu(K,H)=-\sum_{K< T\leq H}\mu(K,T).\]
Equivalently, $\mu$ is defined as the inverse of the zeta function in the incidence algebra of $\SG$. It is easy to see that $\mu(K,H)=-1$ if $K$ is a maximal subgroup of $H$. We write $\mu(G)$ instead of $\mu(1,G)$. Moreover, if $N\unlhd G$, then $\mu(N,G)=\mu(G/N)$ since $\mathcal{S}_{G/N}$ and $\{T\in \SG \mid N\leq T\}$ are isomorphic as posets. The following proposition collects some properties of $\mu$ proved by Hawkes, Isaacs and Özaydin in \cite{HaIsOz}.

\begin{prop}[\cite{HaIsOz}] \thlabel{HaIsOz}
    Let $G$ be a finite group.
    \begin{enumerate}
        \item If $N\unlhd G$, let $\mathcal{K}(N,G)$ denote the set of complements of $N$ in $G$. Then
        $$\mu(G)=\mu(G/N)\left(\sum_{X\in \mathcal{K}(N,G)}\mu(X,G)\right).$$
        \item If $G$ is a $p$-group of order $p^r$, then $\mu(G)\neq 0$ if and only if $G$ is an elementary abelian $p$-group. In this case, $\mu(G)=(-1)^rp^{\binom{r}{2}}$.
        \item If $H\leq G$, then $[N_G(H):H]$ divides $\mu(H,G)[G:HG']_0$.
    \end{enumerate}
\end{prop}

Next, we establish a vanishing criterion for the Möbius function in terms of the indices determined by the subgroups $HG^{[n]}$. This will be useful in the proof of our main result.

\begin{lemma}\thlabel{muvanish}
    Let $H\unlhd G$ and $n$ be a positive integer. If $\left[G:HG^{[n]}\right]$ and $\left[HG^{[n]}:HG'\right]$ are not coprime, then $\mu(H,G)=0$.
\end{lemma}
\begin{proof}
    The quotient $G/HG'$ is abelian and $\left(G/HG'\right)^{[n]}=HG^{[n]}/HG'$. By part 1 of \thref{HaIsOz}, we have $\mu(H,G)=\mu(HG',G)\cdot T$ for some integer $T$. Since $\left[G:HG^{[n]}\right]$ and $\left[HG^{[n]}:HG'\right]$ have a common prime divisor $p$, it follows from \thref{discardprimes} that $G/HG'$ has a $p$-element of order greater than $p$. Consequently, there is a subgroup $HG'\leq M\leq G$ such that $G/M$ is cyclic of order $p^2$. Applying part 1 of \thref{HaIsOz} again, $\mu(HG',G)=\mu(M,G)\cdot U$ for some integer $U$. Finally, $\mu(M,G)=\mu(G/M)=0$ by part 2 of \thref{HaIsOz}, and the result follows.
\end{proof}

\subsection{Fibered Burnside rings}

Given a finite group $G$ and an abelian group $A$, an \textit{$A$-fibered $G$-set} is a $G\times A$-set on which $A$ acts freely and with finitely many orbits. The category of $A$-fibered $G$-sets and $G\times A$-equivariant functions, denoted by $_G\mathrm{set}^A$, is a symmetric semiring category for disjoint unions and a tensor product $\otimes_A$. The \textit{$A$-fibered Burnside ring of $G$} is defined as the Grothendieck ring of the category $_G\mathrm{set}^A$, denoted by $\BAG$. 

The following presentation of $\BAG$ is often more convenient for computations: let $\MAG$ be the set of \textit{$A$-monomial pairs} or \textit{$A$-subcharacters} of $G$, that is, pairs of the form $(H,\phi)$ where $H\leq G$ and $\phi\in \mathrm{Hom}(H,A)$; the group $G$ acts on the set $\MAG$ by conjugation, and the stabilizer and orbit of a pair $(H,\phi)$ are denoted by $N_G(H,\phi)$ and $[H,\phi]_G$, respectively. As the transitive $A$-fibered $G$-sets are in one-to-one correspondence with the orbits of $\MAG$, the ring $\BAG$ is the free abelian group on $\MAG/G$, with multiplication determined on the basis by
$$[H,\phi]_G[K,\psi]_G=\sum_{HgK\in H\backslash G/K}\left[H\cap {^gK},\phi|_{H\cap {^gK}}{^g\psi}|_{H\cap {^gK}}\right]_G,$$
and with $[G,1]_G$ acting as multiplicative identity; $\MAG/G$ is called the \textit{monomial basis} of $\BAG$. It is known that $B^A(G)\cong B^{\mathrm{tor}_{\mathrm{exp}(G)}(A)}(G)$, and moreover, $\BAG$ has finite rank if and only if $\mathrm{tor}_{\mathrm{exp}(G)}(A)$ is finite. Throughout this paper, we assume that $\mathrm{tor}_{\mathrm{exp}(G)}(A)$ is finite. 

For each $H\leq G$, there is a \textit{restriction map} $\mathrm{res} ^G_H:\BAG\longrightarrow \BAH$ defined on the basis by 
$$\mathrm{res}^G_H\left([K,\psi]_G\right)=\sum_{HgK\in H\backslash G/K}\left[H\cap{^gK},{^g\psi|_{H\cap {^gK}}}\right]_H,$$ 
as well as an \textit{induction map} $\mathrm{ind}^G_H:\BAH\longrightarrow \BAG$, $[K,\psi]_H\mapsto [K,\psi]_G$. For $H\leq G$ and $g\in G$, there is also a \textit{$g$-conjugation map} $\mathrm{c}_{g,H}:\BAH \longrightarrow B^A({^gH})$, given by $[K,\psi]_H\mapsto {^g[K,\psi]_H}=\left[{^gK},{^g\psi}\right]_{^gH}$. The collection of the rings $\BAH$ for $H\leq G$ together with induction, restriction and conjugation maps is a Green functor over $G$ in the usual sense.

The ordinary Burnside ring $B(G)$ embeds naturally into $B^A(G)$ via the ring homomorphism $[G/H]\mapsto[H,1]_G$, admitting a retraction $[H,\phi]_G\mapsto [G/H]$. There is also a surjective ring homomorphism $\pi_G:\BAG\longrightarrow \Z\mathrm{Hom}(G,A)$ mapping $[K,\phi]_G$ to $\phi$ if $K=G$, and to $0$ otherwise; this map admits a section $\phi\mapsto [G,\phi]_G$. 

The group $G$ acts on the product of group rings $\prod_{H\leq G}\Z\mathrm{Hom}(H,A)$ through the conjugation maps, and the \textit{ghost ring} of $\BAG$ is defined as the ring of fixed points 
$$\gBAG:=\left(\prod_{H\leq G}\Z\mathrm{Hom}(H,A)\right)^G.$$ 
The product map $\rho=\prod_{H\leq G}\pi_H\circ \mathrm{res}^G_H$ is known as the \textit{mark morphism} and embeds $\BAG$ into its ghost ring. It is explicitly given on the basis by
\[\rho\left([K,\psi]_G\right)= \left(\sum_{\substack{
gK\in G/K,\\
H\leq {^gK}
}}{^g\psi|_H}\right)_{H\leq G}.\]

We write $B^A_R(G)$ for the scalar extension of $\BAG$ to a commutative ring $R$, regarded as an $R$-algebra, and $\widetilde{B^A_R}(G)$ for the corresponding extension of its ghost ring. The extensions of the maps defined in this section are denoted the same. \textit{Boltje's Integrality Criterion} \cite[Theorem 2.2]{Bol1} characterizes the image of the extended mark morphism in $\widetilde{B^A_R}(G)$ as follows. Let
$$y=\left(\sum_{\phi\in \mathrm{Hom}(H,A)}a_{H,\phi}\phi\right)_{H\leq G}$$ 
be an element of $\widetilde{B^A_R}(G)$. Then $y$ belongs to $\rho\left(B^A_R(G)\right)$ if and only if, for every $(K,\psi)\in \MAG$, the congruence
\begin{align}
    \sum_{
    \substack{
    (H,\phi)\in \mathcal{M}^A_{N_G(K,\psi)},\\
    (K,\psi)\leq (H,\phi)
    }}
    \mu(K,H)a_{H,\phi}\; \equiv \; 0 \; \pmod{\left[N_G(K,\psi):K\right]}
\end{align}
holds. Here, $(K,\psi)\leq (H,\phi)$ means that $K\leq H$ and $\psi=\phi|_K$. This defines a partial order on $\MAG$ respecting the action of $G$.

Let $\K$ be a field of characteristic zero that is a splitting field for every $\mathrm{Hom}(H,A)$, with $H\leq G$. As noted in \cite[\S 3]{BolYil}, this is equivalent to requiring that $\K$ contains a root of unity of order $\mathrm{exp}\left(\mathrm{tor}_{\mathrm{exp}(G)}(A)\right)$. We consider the set $\DAG$ consisting of pairs $(H,\Phi)$ where $H\leq G$ and $\Phi \in \mathrm{Hom}(\mathrm{Hom}(H,A),\K^{\times})$; again, $G$ acts on this set by conjugation and the stabilizer and orbit of a pair $(H,\Phi)$ are denoted by $N_G(H,\Phi)$ and $[H,\Phi]_G$, respectively. If $A=C_n$, then $\DnG$ can be identified with the set of pairs $\left(H,hH^{[n]}\right)$, $H\leq G$ and $h\in H$. For each $(H,\Phi)$, the composite $\sGAHPhi:=\Phi\circ \pi_H\circ \mathrm{res}^G_H$ is a ring homomorphism from $\BAG$ to $\K$. By \cite{BolYil}, two of such maps $\sGAHPhi$ and $s^{G,A}_{K,\Psi}$ are equal if and only if $(H,\Phi)=_G (K,\Psi)$. Moreover, these maps account for all elements of $\mathrm{Hom}_{\mathrm{ring}}\left(\BAG,\K\right)$, and $\BAKG$ and $\prod_{[H,\Phi]_G\in \DAG/G}\K$ are isomorphic as $\K$-algebras. 

The primitive idempotents of $\BAKG$ are in bijection with the orbits of $\DAG$. More precisely, for each $(H,\Phi)\in \DAG$ there is a primitive idempotent $\eGAHPhi$ determined by 
\[s^{G,A}_{K,\Psi}\left(\eGAHPhi\right)=\begin{cases}
    1 &\text{if }(K,\Psi)=_G (H,\Phi),\\
    0 &\text{else.}
\end{cases}\]
In terms of the standard basis,
\begin{align}
    \eGAHPhi=\frac{1}{|N_G(H,\Phi)||\mathrm{Hom}(H,A)|}\sum_{\substack{
    K\leq H,\\
    \phi\in \mathrm{Hom}(H,A)
    }
    }|K|\mu(K,H)\Phi(\phi^{-1})[K,\phi|_K]_G,
\end{align}
as showned in \cite{BolYil}. If $\mathcal{O}$ is the ring of algebraic integers of $\K$, the \textit{conductor} of a primitive idempotent $\eGAHPhi$ is the least positive integer $t$ such that $t\eGAHPhi$ belongs to $\BAOG$, denoted by $\cGAHPhi$. In the case of cyclic fiber, we write $\eGnHh$ instead of $e^{G,C_n}_{H,\mathrm{ev}_h}$, and $\cGnHh$ for its conductor. The conductors do not depend on the choice of $\K$, since the coefficients of the primitive idempotents with respect to the monomial basis are algebraic numbers. 

As in \cite{Muller}, the behavior of certain idempotents with respect to induction and restriction will play a role in the proof of the main result. The following lemma collects the required properties. The proof is straightforward and follows the same ideas as the corresponding results in \cite{GarRag} and \cite{Muller}, so we omit it.

\begin{lemma}
    Let $(H,\Phi)\in \mathcal{D}^A_G$. Then:
    \begin{enumerate}
        \item \[\mathrm{res}^G_H\left(e^{G,A}_{H,\Phi}\right)=\sum_{\substack{
        \Psi\in \mathrm{Hom}\left(\mathrm{Hom}(H,A), \K^{\times}\right),\\
        [H,\Psi]_G=[H,\Phi]_G
        }}e^{H,A}_{H,\Psi},\]
        \item \[\mathrm{ind}^G_H\mathrm{res}^G_H\left(e^{G,A}_{H,\Phi}\right)=[N_G(H):H]e^{G,A}_{H,\Phi},\]
        \item \[\mathrm{ind}^G_H\left(e^{H,A}_{H,\Phi}\right)=\left[N_G(H,\Phi):H\right]e^{G,A}_{H,\Phi}.\]
    \end{enumerate}
\end{lemma}

\section{The conductor of a primitive idempotent}

In the following, $A$ is an abelian group such that $\mathrm{tor}_{\mathrm{exp}(G)}(A)$ is finite. Let $\K$ be a field of characteristic zero that is a splitting field for every $\mathrm{Hom}(H,A)$, with $H\leq G$, and let $\mathcal{O}$ be its ring of algebraic integers. Consider the primitive idempotent of $B_{\Q}(G)$ corresponding to the $G$-th mark,
$$e^{B(G)}_G=\sum_{H\in [\SG]}a_H[G/H],$$
with $a_H\in \Q$. The following characterization of a primitive idempotent $\eGAGPhi$ is a generalization of that presented by Müller in \cite[\S 3]{Muller}.

\begin{lemma}\thlabel{newidempotent}
    Let $\Phi \in \mathrm{Hom}\left(\mathrm{Hom}(G,A),\K^\times\right)$. Then
    \begin{align}
        \eGAGPhi=\frac{1}{|\mathrm{Hom}(G,A)|}\sum_{\lambda\in \mathrm{Hom}(G,A)}\Phi\left(\lambda^{-1}\right)x_{\lambda},
    \end{align}
    where 
    $$x_{\lambda}:=\sum_{H\in [\SG]}a_H[H,\lambda|_H]_G.$$
\end{lemma}
\begin{proof}
    It is straightforward that $\pi_K\left([H,\lambda|_H]_G\right)=\phi_K(G/H)\lambda|_K$ for $K\leq G$, hence
    \begin{align*}
        \pi_K\circ \mathrm{res}^G_K\left(x_{\lambda}\right)=\phi_K\left(e^{B(G)}_G\right)\lambda|_K=\begin{cases}
            \lambda &\text{if }K=G,\\
            0 &\text{otherwise.}
        \end{cases}
    \end{align*}
    For arbitrary $(K,\Psi)\in \DAG$, we have the following: if $K$ is proper, then
    $$s^{G,A}_{K,\Psi}\left(\sum_{\lambda\in \mathrm{Hom}(G,A)}\Phi\left(\lambda^{-1}\right)x_{\lambda}\right)=\sum_{\lambda\in \mathrm{Hom}(G,A)}\Phi\left(\lambda^{-1}\right)s^{G,A}_{K,\Psi}\left(x_{\lambda}\right)=0;$$
    if $K=G$, then
    \begin{align*}
        s^{G,A}_{G,\Psi}\left(\sum_{\lambda\in \mathrm{Hom}(G,A)}\Phi\left(\lambda^{-1}\right)x_{\lambda}\right) & =\sum_{\lambda\in \mathrm{Hom}(G,A)}\Phi\left(\lambda^{-1}\right)\Psi\left(\lambda\right)\\
        &=\begin{cases}
            |\mathrm{Hom}(G,A)| &\text{if }\Psi=\Phi,\\
            0 &\text{otherwise}.
        \end{cases}
    \end{align*}
    by row orthogonality relations for the group $\mathrm{Hom}(G,A)$, proving that
    \[|\mathrm{Hom}(G,A)|\eGAGPhi=\sum_{\lambda\in \mathrm{Hom}(G,A)}\Phi\left(\lambda^{-1}\right)x_{\lambda}.\]
\end{proof}

For a positive integer $n$, we set
$$r^n(G):=\frac{\left[G^{[n]}:G'\right]_0}{\gcd\left(n,\left[G^{[n]}:G'\right]\right)_0},$$
that is, $r^n(G)$ is the product of the prime divisors of $\left[G^{[n]}:G'\right]$ which are coprime to $n$. In general, we set 
$$r^A(G):=r^{\mathrm{exp}\left(\mathrm{tor}_{\mathrm{exp}(G)}(A)\right)}(G).$$
This notation is useful, since the exponent of $\mathrm{tor}_{\mathrm{exp}(G)}(A)$ may change as $G$ varies.

\begin{theorem}
    Let $G$ be a finite group, and let $A$ be an abelian group such that
    $\mathrm{tor}_{\mathrm{exp}(G)}(A)$ is finite. Let $\K$ be a field of characteristic zero that is a splitting field for every $\mathrm{Hom}(H,A)$, with $H\leq G$. For $(H,\Phi)\in\mathcal{D}^A_G$, the conductor of the primitive idempotent $e^{G,A}_{H,\Phi}$ of $B^A_{\K}(G)$ is
    \begin{align}
        c^{G,A}_{H,\Phi}=\left[N_G(H,\Phi):H\right]\left|\mathrm{Hom}(H,A)\right|r^A(H).
    \end{align}
\end{theorem}
\begin{proof}
    Let
    $$\mathrm{tor}_{\mathrm{exp}(G)}(A)=C_{d_1}\times C_{d_2}\times \cdots \times C_{d_t}$$
    be a decomposition of $\mathrm{tor}_{\mathrm{exp}(G)}(A)$ into invariant factors, where $d_i\mid d_{i+1}$ for $i=1,...,t-1$, and $d_t =\mathrm{exp}\left(\mathrm{tor}_{\mathrm{exp}(G)}(A)\right)$. We first consider the case $H=G$, and set
    \[y :=\rho\left(\left|\mathrm{Hom}(G,A)\right|\eGAGPhi\right)=\left(\sum_{\phi\in \mathrm{Hom}(K,A)}a_{K,\phi}\phi\right)_{K\leq G},\]
    which lies in $\gBAOG$ by \thref{newidempotent}. We test Boltje's integrality criterion on $y$. Because the $K$-th entry of $y$ is zero for any proper subgroup $K$, we only need to check the congruences for $(K,\phi)$ such that $N_G(K,\phi)=G$ and $\phi\in \mathrm{Im}\left(\mathrm{res}^G_K\right)$, that is, $K\unlhd G$ and $\phi=\lambda|_K$ for some $\lambda \in \mathrm{Hom}(G,A)$. For such a pair, we have
    \begin{align*}
        \sum_{
    \substack{
    (L,\psi)\in \MAG,\\
    (K,\phi)\leq (L,\psi)
    }}
    &\mu(K,L)a_{L,\psi} =\sum_{\lambda\in \left(\text{res}^G_K\right)^{-1}(\phi)}\mu(K,G)\Phi(\lambda^{-1})\\
        &=\begin{cases}
            \mu(K,G){|\mathrm{ker}\left(\mathrm{res}^G_K\right)|}\Phi(\lambda^{-1}_0) &\text{if }\text{ker}\left(\text{res}^G_K\right)\leq \text{ker}(\Phi),\\
            0 &\text{otherwise,}
        \end{cases}
    \end{align*}
    for some $\lambda_0\in \left(\mathrm{res}^G_K\right)^{-1}\left(\phi\right)$. Since there is a natural identification $\mathrm{Hom}(L,A) \cong \mathrm{Hom}\left(L,\mathrm{tor}_{\mathrm{exp}(G)}(A)\right)$ for $L\leq G$, we have $\left|\mathrm{ker}\left(\mathrm{res}^G_K\right)\right|=\prod_{i=1}^t\left[G:KG^{[d_i]}\right]$ by \thref{MullerLemma}, while $\Phi(\lambda_0^{-1})$ is a root of unity. Moreover, every $\left[G:KG^{[d_i]}\right]$ is a divisor of $\left[G:KG^{[d_t]}\right]$. If $\left[G:KG^{[d_t]}\right]$ and $\left[KG^{[d_t]}:KG'\right]$ are not coprime, then $\mu(K,G)=0$ by \thref{muvanish}. On the other hand, if $\gcd\left(\left[G:KG^{[d_t]}\right],\left[KG^{[d_t]}:KG'\right]\right)=1$, every prime divisor of $\left[KG^{[d_t]}:KG'\right]$ is a divisor of $\left[G^{[d_t]}:G'\right]$ that does not divide $d_t$ by \thref{discardprimes}, so it is a divisor of $r^{d_t}(G)=r^A(G)$. Hence, $\mu(K,G)\left|\mathrm{Hom}(G,A)\right|r^A(G)$ is multiple of $\mu(K,G)\left[G:KG'\right]_0$, and therefore a multiple of $[N_G(K):K]=[G:K]$, by part 3 of \thref{HaIsOz}. Since $N_G(K,\phi)=G$, Boltje's integrality criterion implies that
    $$r^A(G)y=\rho\left(\left|\mathrm{Hom}(G,A)\right|r^A(G)\eGAGPhi\right) \in \rho\left(\BAOG\right),$$
    hence $\left|\mathrm{Hom}(G,A)\right|r^A(G)\eGAGPhi \in \BAOG$. For general $(H,\Phi)$, we have that
    \begin{align*}
        \left[N_G(H,\Phi):H\right]\left|\mathrm{Hom}(H,A)\right| r^A(H) e^{G,A}_{H,\Phi} = \mathrm{ind}^G_H\left(\left|\mathrm{Hom}(H,A)\right|r^A(H)e^{H,A}_{H,\Phi}\right) 
    \end{align*}
    belongs to $\BAOG$ since $\left|\mathrm{Hom}(H,A)\right|r^A(H)e^{H,A}_{H,\Phi}$ belongs to $B^A_{\mathcal{O}}(H)$. This shows that $\left[N_G(H,\Phi):H\right]\left|\mathrm{Hom}(H,A)\right|r^A(H)$ divides $\cGAHPhi$.

    Conversely, note that $\mathrm{tor}_{\mathrm{exp}(H)}(A)$ is also finite, and consider a decomposition
    $$\mathrm{tor}_{\mathrm{exp}(H)}(A)=C_{f_1} \times C_{f_2}\times \cdots \times C_{f_s}$$ 
    into invariant factors, where $f_i\mid f_{i+1}$ for $i=1,...,s-1$ and $f_s=\mathrm{exp}\left(\mathrm{tor}_{\mathrm{exp}(H)}(A)\right)$. Note that the coefficient of $[H,1]_G$ in $\eGAHPhi$ with respect to the monomial basis is $\left([N_G(H,\Phi):H]|\mathrm{Hom}(H,A)|\right)^{-1}$, therefore $[N_G(H,\Phi):H]|\mathrm{Hom}(H,A)|$ divides $\cGAHPhi$. We now prove that $r^A(H)=r^{f_s}(H)$ divides 
    $$u(H,\Phi):=\frac{\cGAHPhi}{\left[N_G(H,\Phi):H\right]\left|\mathrm{Hom}(H,A)\right|}.$$ 
    Let $p$ be a prime divisor of $\left[H^{[f_s]}:H'\right]$ that does not divide $f_s$. For any $K\leq H$ such that $[H:K]=p$, we have $\left[H:KH^{[f_s]}\right]=1$, hence $\left[H:KH^{[f_i]}\right]=1$ for each $i$, and so
    $$\left|\mathrm{ker}\left(\mathrm{res}^H_K\right)\right|=\prod_{i=1}^s\left[H:KH^{[f_i]}\right]=1,$$
    proving that $\mathrm{res}^H_K:\mathrm{Hom}(H,A)\to \mathrm{Hom}(K,A)$ is injective. Let 
    $$\mathcal{I}_p=\left\{K\leq H \mid [H:K]=p\right\}.$$ 
    We have $\left|\mathcal{I}_p\right|\equiv 1 \pmod{p}$, and $\mathcal{I}_p$ splits into equivalence classes under the relation $=_G$. Therefore, there is some $K\in \mathcal{I}_p$ such that $p$ does not divide 
    $$m=\# \{L\in \mathcal{I}_p \mid L=_GK\}.$$
    For such $K$, the coefficient of $[K,1]_G$ in $\left[N_G(H,\Phi):H\right]\left|\mathrm{Hom}(H,A)\right|\eGAHPhi$ with respect to the monomial basis is
    $$\frac{1}{|H|}\sum_{\substack{
    L\leq H,\\
    L=_GK,\\
    \phi\in \mathrm{ker}\left(\mathrm{res}^H_L\right)
    }
    }|L|\mu(L,H)\Phi\left(\phi^{-1}\right)=-\frac{|K|}{|H|}\sum_{\substack{
        L\leq H,\\
        L=_GK}
        }\Phi(1)=-\frac{m}{p}.$$
    Since $u(H,\Phi)\left(-\frac{m}{p}\right)$ must be an integer, we conclude that $p$ must divide $u(H,\Phi)$, completing the proof.
\end{proof}

\begin{cor}
    Let $G$ be a finite group, and let $n$ be a positive integer. For $\left(H,hH^{[n]}\right)\in \DnG$, the conductor of $\eGnHh$ is
    \begin{align}
        \cGnHh=\left[N_G\left(H,hH^{[n]}\right):H^{[n]}\right]r^n(H).
    \end{align}
\end{cor}

For both the Burnside ring and the monomial representation ring, the order of a group can be easily recovered from their conductors \cite[Proposition 3.7]{Muller}. This also holds for all fibered Burnside rings with cyclic fiber group. 

\begin{prop}
    Let $n$ be a positive integer. For a finite group $G$, we have
    \begin{align*}
        |G| &=\max\left\{\cGnHh\mid \left(H, hH^{[n]}\right)\in \DnG\right\}\\
        &=\min \left\{t\in \Z_{>0} \mid t \eGnHh\in \BAOG, \forall \left(H,hH^{[n]}\right)\in \DnG\right\}.
    \end{align*}
\end{prop}
\begin{proof}
Take a characteristic zero field $\K$ containing a primitive $n$-th root of unity. The conductors of the primitive idempotents of $B^{C_n}_{\K}(G)$ do not depend on the field $\K$, and $\left[N_G\left(H,hH^{[n]}\right):H^{[n]}\right]r^n(H)$ divides $|G|$ for every pair $\left(H,hH^{[n]}\right)$. As $\cGnHh=|G|$ is reached at least once, this proves both equalities.
\end{proof}

The previous argument fails for more general fibers, as shown in the following example.

\begin{exx}
    Let $G=C_4$, $H\leq C_4$ be the subgroup of order $2$ and $A=C_2\times C_2$. Then $N_G(H,1)=N_G(H)=C_4$, $\left|\mathrm{Hom}(H,C_2\times C_2)\right|=4$, and $r^{C_2\times C_2}(H)=1$. Thus, we have
    $$c^{C_4,C_2 \times C_2}_{H,1}=\left[C_4:H\right]\left|\mathrm{Hom}(H,C_2\times C_2)\right|r^{C_2\times C_2} (H)=8,$$
    which exceeds $|C_4|=4$.
\end{exx}

The idempotents of $B(G)$ whose conductors are exactly $|G|$ play a central role in the proof of the existence of normalizing isomorphisms, with important consequences on the isomorphism problem of the Burnside ring \cite{RagVal}. In the case of $D(G)$, we see that $G$ is abelian if and only if all conductors are equal to $|G|$ \cite{Muller}. More generally, we can characterize pairs $\left(H,hH^{[n]}\right)$ for which $\cGnHh=|G|$. 

\begin{prop}
    Let $n$ be a positive integer. For $\left(H,hH^{[n]}\right)\in \DnG$, $\cGnHh=|G|$ if and only if $H$ is a normal abelian subgroup such that $H^{[n]}$ has square-free order coprime to $n$ and $hH^{[n]}\in Z\left(G/H^{[n]}\right)$.
\end{prop}
\begin{proof}
    Note that $\left[N_G\left(H,hH^{[n]}\right):H^{[n]}\right]r^n(H)=[G:1]$ implies $N_G\left(H,hH^{[n]}\right)=G$, and so $N_G(H)=G$ and $hH^{[n]} \in Z\left(G/H^{[n]}\right)$. Then $\left[G:H^{[n]}\right]r^n(H)=[G:1]$ implies $\left[H^{[n]}:1\right]=r^n(H)$, the latter is a divisor of $\left[H^{[n]}:H'\right]$, thus implying $H'=1$ and proving the assertion about the order of $H^{[n]}$. Conversely, suppose that $H$ is normal and abelian, and that $H^{[n]}$ has square-free order coprime to $n$. For $hH^{[n]}\in Z\left(G/H^{[n]}\right)$, we have $N_G\left(H,hH^{[n]}\right)=G$, and, since $H'=1$, we have $r^n(H)=\left|H^{[n]}\right|$. Consequently, $c^{G,n}_{H,h}=\left[G:H^{[n]}\right]\left|H^{[n]}\right|=|G|$.
\end{proof}

Specializing to $p$-groups, we get the following.

\begin{cor}
    For a prime number $p$, let $P$ be a $p$-group and let $q$ be a power of $p$. For $\left(H,hH^{[q]}\right)\in \mathcal{D}^{C_q}_P$, $c^{P,q}_{H,h}=|P|$ if and only if $H$ is a normal abelian subgroup such that $H^{[q]}=1$ and $h\in Z(P)$.
\end{cor}
\begin{proof}
    This follows immediately from the previous result, since in this case $r^q(H)=1$.
\end{proof}

\end{document}